\let\OToneAccents=\relax
\magnification=\magstep1
\frenchspacing
\baselineskip=16truept
\font\douze=cmr10 at 12pt

\font\grande=cmb10 at 16truept
\def\titre#1{{\OToneAccents\noindent\grande #1}}

\def \chapter#1{\vfill\eject\ifodd\pageno\else\ \vfill\eject\fi\centerline{\grande #1}\bigskip}

\font\tendo=wncyr10 at 12pt
\font\sevendo=wncyr7
\newfam\dofam 
\textfont\dofam = \tendo 
\scriptfont\dofam= \sevendo

\font\tendb=msbm10 
\font\sevendb=msbm7
\newfam\dbfam 
\textfont\dbfam = \tendb 
\scriptfont\dbfam= \sevendb
\def\db {\fam\dbfam\tendb}
\font\tenrsfs=rsfs10 
\font\sevenrsfs=rsfs7
\newfam\scrfam 
\textfont\scrfam = \tenrsfs
\scriptfont\scrfam= \sevenrsfs

\def\R{{\db R}}
\def\Q{{\db Q}}

\def\Z{{\db Z}}

\def\Cl{{\db Cl}}

\def\Cl{\mathop{\rm Cl}\nolimits}

\font\tendb=msbm10 
\font\sevendb=msbm7
\newfam\dbfam 
\textfont\dbfam = \tendb 
\scriptfont\dbfam= \sevendb
\def\db {\fam\dbfam\tendb}
\font\tenrsfs=rsfs10 
\font\sevenrsfs=rsfs7
\newfam\scrfam 
\textfont\scrfam = \tenrsfs
\scriptfont\scrfam= \sevenrsfs

\def\R{{\db R}}
\def\Q{{\db Q}}
\def\Z{{\db Z}}

\def\FF{{\db F}}

\def\mod{\hbox{\rm \ mod.\ }}

\def\deg{\mathop{\rm deg}\nolimits}

\font\teufm=eufm10
\font\seufm=eufm10 at 7pt
\font\sseufm=eufm10 at 6pt
\newfam\fameufm
\textfont\fameufm=\teufm
\scriptfont\fameufm=\seufm
\scriptscriptfont\fameufm=\sseufm
\def\goth{\fam\fameufm\teufm}

\def\wp{{\goth p}}

\def\qq{{\goth q}}

\def\NN{\mathop{\rm N}\nolimits}

\def\mod{\hbox{\rm \ mod.\ }}

\def\Cl{\mathop{\rm Cl}\nolimits}

\bigskip\bigskip

\bigskip\bigskip

\def\og{\leavevmode\raise.3ex\hbox{$\scriptscriptstyle\langle\!\langle\,$}}
\def\fg{\leavevmode\raise.3ex\hbox{$\scriptscriptstyle\,\rangle\!\rangle$}}

\centerline {\titre{Contre-exemples au principe de Hasse}}
\medskip
\centerline {\titre{pour les courbes de Fermat }}
\bigskip 
\centerline {Alain Kraus}
\bigskip 
\medskip
{\bf{Abstract.}} Let $p$ be an odd  prime number. In this paper, we are concerned with the behaviour of Fermat curves defined over $\Q$ given by equations
$ax^p+by^p+cz^p=0$, with respect to the local-global Hasse principle. It is conjectured that there exist infinitely many Fermat curves of exponent $p$ 
which are counterexamples to the Hasse principle. It is a consequence of the abc-conjecture if $p\geq 5$. Using  a cyclotomic approach due to H. Cohen and Chebotarev's density theorem, we  obtain a partial result towards this conjecture,  by proving it for $p\leq 19$.
\bigskip

\medskip

\bigskip
\smallskip

\centerline {\douze{Introduction}}
\medskip
Soit $p\geq 3$ un nombre premier. Une courbe de Fermat  $C/\Q$ d'exposant $p$ est d\'efinie par une \'equation de la forme 
$$ax^p+by^p+cz^p=0,$$
o\`u $a,b,c$ sont des entiers rationnels non nuls. 
On s'int\'eresse  dans cet article au comportement de ces  courbes   vis-\`a-vis du principe de Hasse.   
Adoptons la terminologie en vigueur  selon laquelle $C$ pr\'esente une obstruction locale en un nombre premier $\ell$ si $C(\Q_{\ell})$ est vide. On dit   que la courbe $C$ contredit le principe de Hasse  si elle  ne pr\'esente aucune obstruction locale  et si  $C(\Q)$ est vide. La place \`a l'infini n'intervient pas ici car $C(\R)$ est non vide.
\smallskip
 Si $a,b,c$   ne v\'erifient pas de relation lin\'eaire non triviale \`a coefficients dans $\big\lbrace -1,0,1\big\rbrace$, il est tr\`es fr\'equent que $C$ poss\`ede au moins une obstruction locale ([H-K]). On se pr\'eoccupe  ici du probl\`eme d'expliciter des courbes de Fermat contredisant le  principe de Hasse. \`A ma connaissance, l'\'etat de ce probl\`eme est le suivant. Historiquement, le premier exemple qui a \'et\'e d\'ecouvert \`a ce sujet est la courbe d'\'equation 
 $$3x^3+4y^3+5z^3=0,$$ 
explicit\'ee  par   E. S. Selmer en 1951 ([Se]). 
Par ailleurs,
il figure   dans  [H-K] des exemples de courbes de Fermat contredisant ce principe pour tout $p$ au moins $5$   plus petit que $100$. Ils ont \'et\'e obtenus par des techniques modulaires li\'ees aux repr\'esentations galoisennes des points de torsion des courbes elliptiques. H. Cohen a  \'egalement obtenu de tels exemples pour $p\leq 11$ par une approche cyclotomique ([Co-2], cor. 6.4.11). Cela \'etant, on  ne sait pas d\'emontrer que pour tout $p$, il existe une courbe de Fermat d'exposant $p$ contredisant  le principe de Hasse. 
N\'eanmoins, certains r\'esultats   \'etablis dans [H-K] rendent plausible la conjecture sui\-vante~:
\medskip

\proclaim Conjecture. Pour tout nombre premier $p\geq 3$, il existe une infinit\'e de courbes de Fermat d'exposant  $p$,  deux \`a deux non $\Q$-isomorphes,  
contredisant le principe de Hasse.
\medskip

C'est une cons\'equence de la conjecture abc pour  $p\geq 5$ ({\it{loc. cit.}}, prop. 5.1).  
Le nombre premier $p$ \'etant donn\'e, on  obtient  ici un crit\`ere impliquant  cet \'enonc\'e  pour l'exposant $p$. Cela permet d'en d\'eduire  cette conjecture pour $p\leq 19$.
  La m\'ethode que l'on utilise repose sur   l'approche cyclotomique pr\'esent\'ee par H. Cohen dans [Co-2], ainsi que sur  le th\'eor\`eme de densit\'e de Chebotarev. Elle  permet par ailleurs,  pour $p\leq 19$, d'expliciter de nombreuses courbes de Fermat d'exposant $p$ contredisant le principe de Hasse.
\smallskip
Tous les calculs num\'eriques que  ce travail a n\'ecessit\'es ont \'et\'e effectu\'es avec le logiciel de calcul  Pari ([Pa]). 
\smallskip
Je remercie  D.  Bernardi pour les remarques dont il m'a fait part  concernant cet article.
\bigskip

{\douze I. \'Enonc\'e des r\'esultats}
\medskip
Soient $p\geq 3$ un nombre premier et $c\geq 2$ un entier sans puissances $p$-i\`emes.
Pour tout nombre premier $\ell$, distinct de $c$, d\'esignons par    $C_{\ell}/\Q$ la courbe d'\'equation
$$x^p+\ell y^p+cz^p=0.$$
\item{}  Posons $K=\Q\bigl(\root p \of {c}\bigr)$ et notons : 
\medskip
\vskip0pt\noindent
\item{.} $O_K$ son anneau d'entiers,
\smallskip
\vskip0pt\noindent
\item{.}  $f$ l'indice de $\Z\big[\root p \of {c}\big]$   dans $O_K$,
\smallskip
\vskip0pt\noindent
\item{.} $\Cl_K$ le groupe des classes de $K$, 
\smallskip
\vskip0pt\noindent
\item{.} $h_K$ le nombre de classes de $K$.
\medskip
\vskip0pt\noindent
\item{} On supposera  dans toute la suite que   $p$ divise $h_K$.  Notons par ailleurs :
\medskip
\vskip0pt\noindent
\item{.} $e$ l'exposant du groupe ab\'elien $\Cl_K$,
\smallskip
\vskip0pt\noindent
\item{.} $r$ la valuation $p$-adique de $e$,
\smallskip
\vskip0pt\noindent
\item{.} $N$ le nombre de copies de $\Z/p^r\Z$ intervenant dans la d\'ecomposition primaire de $\Cl_K$, 
\smallskip
\vskip0pt\noindent
\item{.}  $S$  l'ensemble des nombres premiers $\ell\not\equiv 1\mod p$, ne divisant pas $f$, tels que la courbe $C_{\ell}/\Q$  
ne pr\'esente aucune obstruction locale. 
\vskip0pt\noindent
\item{} Pour tout $\ell\in S$, il existe un unique id\'eal premier de $O_K$ au-dessus de $\ell$ de degr\'e
r\'esiduel 1 (lemme 1).
\smallskip
\vskip0pt\noindent
\item{.} $S_0$ le sous-ensemble de $S$ form\'e des nombres premiers $\ell$  tels que    la condition suivante soit satisfaite  : 

\item{} soit $\qq$ l'id\'eal premier   de $O_K$ au-dessus de $\ell$  de degr\'e r\'esiduel $1$.  L'id\'eal $\qq^{{e\over p}}$  n'est pas principal.
\vfill\eject

\proclaim Th\'eor\`eme. Supposons que les deux conditions suivantes  soient remplies :
\smallskip
\vskip0pt\noindent
\item{1)}  on a $c^{p-1}\not\equiv 1 \mod p^2$.
\smallskip
\vskip0pt\noindent
\item{2)}  $p$ divise $h_K$.
\smallskip
\vskip0pt\noindent
Alors, pour  tout $\ell\in S_0$  la courbe $C_{\ell}/\Q$ est  un contre-exemple au principe de Hasse.  
Si $S$ a une densit\'e strictement plus grande  que ${1\over p^N}$, l'ensemble $S_0$ est infini, auquel cas la conjecture est vraie pour l'exposant $p$.
 \bigskip

\proclaim Corollaire. Supposons que  $(p,c)$ soit l'un des couples suivants :
$$(3,921),\quad (5,19),\quad (7,13),\quad (11,373),\quad (13,103), \quad (17,1087),\quad (19,37).$$
Alors, $S_0$ est infini. En particulier, la conjecture est vraie pour $p\leq 19$.
\bigskip

{\bf{Remarque 1.}} L'\'enonc\'e du th\'eor\`eme permet d'obtenir de nombreux contre-exem\-ples au principe de Hasse. \`A titre indicatif, pour $(p,c)=(5,19)$, on a $h_K=5$ et  il y a 72 classes de nombres premiers $\ell$ modulo $275$  qui sont dans $S$. Pour chaque nombre premier $\ell$ dans l'une de ces classes, si $\ell$ est dans $S_0$, alors  $C_{\ell}$ est donc un contre-exemple au principe de Hasse. Par exemple, l'ensemble des nombres premiers congrus \`a $7$ modulo $275$ est l'une de ces classes. Il y a 48  nombres premiers plus petits que $10^5$ congrus \`a $7$ modulo $275$, et  il y en a 31   qui  sont dans  $S_0$.  Le plus petit d'entre eux est $1657$ et le plus grand est $95707$.
\medskip
{\bf{Remarque 2.}} Si l'on essaye  d'\'etendre l'\'enonc\'e du corollaire avec  des   nombres premiers $p\geq 23$,  cela devient, a priori,  nettement plus  couteux num\'eriquement.  Par exemple, afin  de pouvoir conclure pour $p=23$, 
il semble qu'il faille  disposer de valeurs de $c$ pour lesquelles $\Cl_K$ contienne un sous-groupe isomorphe \`a $\Z/23\Z\times \Z/23\Z$.
\bigskip

{\douze II. D\'emonstration du th\'eor\`eme}
\medskip

  Commenons par \'etablir trois  lemmes pr\'eliminaires.
  \bigskip

\proclaim Lemme 1. Soit $\ell$ un nombre premier ne divisant pas $f$ tel que  $\ell\not\equiv 1 \mod p$. Il existe un unique id\'eal premier 
de $O_K$ au-dessus de $\ell$ de degr\'e r\'esiduel $1$.

D\'emonstration : Posons $F=X^p-c$. La condition $\ell\not\equiv 1 \mod p$ entra\^\i ne que $1$ est la seule racine $p$-i\`eme de l'unit\'e dans $\FF_{\ell}$, donc le morphisme $\FF_{\ell}^*\to \FF_{\ell}^*$ qui \`a $x$ associe $x^p$ est bijectif.
Par suite, $F$ a une unique racine dans $\FF_{\ell}$. Parce que $\ell$ ne divise pas $f$, cela entra\^\i ne le r\'esultat ([Co-1], th. 4.8.13).
\bigskip

\proclaim Lemme 2. Soit $\ell$ un nombre premier ne divisant pas $cf$ tel que $\ell\equiv 1 \mod p$. Alors, $\ell$ est inerte ou totalement d\'ecompos\'e dans $K$. 

D\'emonstration :  Les racines $p$-i\`emes de l'unit\'e sont dans $\FF_{\ell}$ car $p$ divise $\ell-1$. Si le polyn\^ome $F=X^p-c$ n'a pas de racines dans $\FF_{\ell}$, alors $F$ est irr\'eductible modulo $\ell$ et $\ell$ est donc inerte dans $K$. Sinon, $F$ a toutes ses racines dans $\FF_{\ell}$, donc poss\`ede $p$ racines distinctes dans $\FF_{\ell}$ ($\ell$ ne divise pas $c$) et $\ell$ est  totalement d\'ecompos\'e dans $K$.
\bigskip

\proclaim Lemme 3. La densit\'e des nombres premiers totalement d\'ecompos\'es dans $K$ est ${1\over p(p-1)}$. 

D\'emonstration :  Un nombre premier  est totalement d\'ecompos\'e dans $K$ si et seulement si il est totalement d\'ecompos\'e dans la cl\^oture galoisienne de $K$, 
qui est de degr\'e $p(p-1)$ sur $\Q$, d'o\`u l'assertion ([Ne], chap. V, cor (6.5)).
\bigskip

1) Consid\'erons  alors un nombre premier $\ell\in S_0$. V\'erifions que   $C_{\ell}(\Q)$ est vide, ce qui \'etablira  que $C_{\ell}/\Q$ est un contre-exemple au principe de Hasse.  On utilise pour cela le th\'eor\`eme 6.4.8 de [Co-2].  Ses deux premi\`eres conditions  sont par hypoth\`ese satisfaites. De plus, $\ell$ ne divise pas $f$. Il en r\'esulte que le seul diviseur convenable de $\ell O_K$, au sens de la d\'efinition 6.4.6 de {\it{loc. cit.}},  est l'id\'eal premier $\qq$ 
de $O_K$ au-dessus de $\ell$ de degr\'e r\'esiduel $1$ (lemme 1). L'id\'eal $\qq^{{e\over p}}$ n'\'etant pas  principal, et $\ell$ \'etant distinct de $c$, cela entra\^\i ne  notre assertion.
\bigskip

2) Supposons $S$ de densit\'e strictement plus grande que ${1\over p^N}$. D\'emontrons  que $S_0$ est infini et par suite que la conjecture est vraie pour l'exposant $p$.
\bigskip

\proclaim Proposition 1. Soit $A$ un  ensemble de nombres premiers non congrus \`a $1$ modulo $p$.  Soit $B$ 
l'ensemble des  id\'eaux premiers de $O_K$ de degr\'e r\'esiduel $1$ au-dessus d'un nombre premier de $A$. Alors, si  $A$ a une densit\'e, il en est de m\^eme de $B$ et elles sont \'egales. 

D\'emonstration :  Pour tout id\'eal premier $\wp$ de $O_K$,  notons $\NN \wp$ sa norme sur $\Q$ et $\deg \wp$ son degr\'e r\'esiduel. Si $\wp$ est au-dessus du nombre premier $\ell$, on a $\NN \wp=\ell^{\deg \wp}$. Supposons $A$ de  densit\'e $d$. 
Pour tout $x>0$ posons 
$$u(x)={\Big| \big\lbrace \wp\in B \ | \ \NN \wp\leq x\big\rbrace \Big|\over \Big| \big\lbrace \wp \ | \ \NN \wp\leq x\big\rbrace \Big|}.$$
Il s'agit d'\'etablir que $u(x)$ a une limite \'egale \`a $d$ quand $x$ tend vers l'infini.
On a 
$$u(x)={\Big| \big\lbrace \wp\in B \ | \ \NN \wp\leq x\big\rbrace \Big|\over \Big| \big\lbrace \wp \ | \ \NN \wp\leq x, \ \deg \wp=1\big\rbrace \Big|}.{\Big| \big\lbrace \wp \ | \ \NN \wp\leq x,\  \deg \wp=1\big\rbrace \Big|\over \Big| \big\lbrace \wp \ | \ \NN \wp\leq x\big\rbrace \Big|}.$$
D'apr\`es les lemmes 1 et 2, il existe une constante   $c_1$ telle  que l'on ait 
$$\Big| \big\lbrace \wp \ | \ \NN \wp\leq x, \ \deg \wp=1\big\rbrace \Big|=\Big| \big\lbrace \ell \ | \ \ell\leq x, \  \ell\not\equiv 1\mod p\big\rbrace \Big|$$
$$+\ p.\Big| \big\lbrace \ell \ | \ \ell\leq x, \ \ell\ \hbox{totalement d\'ecompos\'e dans}\ K\big\rbrace \Big|+c_1.$$
Notons $\pi(x)$ le nombre des nombres premiers plus petits que $x$ et posons 
$$v(x)={\Big| \big\lbrace \wp \ | \ \NN \wp\leq x, \ \deg \wp=1\big\rbrace \Big|\over \pi(x)}.$$
D'apr\`es le th\'eor\`eme de Dirichlet et le lemme 3, la fonction $v(x)$ a une limite quand $x$ tend vers l'infini, qui vaut  ainsi
$${p-2\over p-1}+p.{1\over p(p-1)}=1.$$
On \'ecrit alors que l'on a
$${\Big| \big\lbrace \wp\in B \ | \ \NN \wp\leq x\big\rbrace \Big|\over \Big| \big\lbrace \wp \ | \ \NN \wp\leq x, \ \deg \wp=1\big\rbrace \Big|}={\Big| \big\lbrace \wp\in B \ | \ \NN \wp\leq x\big\rbrace \Big|\over \pi(x)}.{1\over v(x)}.$$
Les nombres premiers de $A$ \'etant non congrus \`a 1 modulo $p$, il  existe  une constante $c_2$  telle que l'on ait (lemme 1)
$$\Big| \big\lbrace \wp\in B \ | \ \NN \wp\leq x\big\rbrace \Big|=\Big| \big\lbrace \ell \ | \ \ell\leq x, \  \ell\in S\big\rbrace \Big|+c_2.$$
Puisque $A$ est de densit\'e $d$, on obtient   
$$\lim_{x\to \infty}{\Big| \big\lbrace \wp\in B \ | \ \NN \wp\leq x\big\rbrace \Big|\over \pi(x)}=d.$$
Par suite, on a
$$\lim_{x\to \infty} {\Big| \big\lbrace \wp\in B \ | \ \NN \wp\leq x\big\rbrace \Big|\over \Big| \big\lbrace \wp \ | \ \NN \wp\leq x, \ \deg \wp=1\big\rbrace \Big|}=d.$$
L'ensemble des id\'eaux premiers de $O_K$ de degr\'e 1 \'etant de densit\'e $1$, la limite de $u(x)$ vaut donc $d$,
d'o\`u le r\'esultat.
\bigskip
{\bf{Remarque 3.}}  L'\'enonc\'e de la proposition 1  est faux si on enl\`eve l'hypoth\`ese $\og$non congrus \`a $1$ modulo $p$$\fg$, comme on le constate  en prenant pour $A$ l'ensemble des nombres premiers inertes dans $K$. En effet, d'apr\`es les lemmes 1 et 2 ainsi que  le th\'eor\`eme  de Dirichlet, la densit\'e des nombres premiers totalement d\'ecompos\'es dans $K$ ou inertes dans $K$ est ${1\over p-1}$. La densit\'e des nombres premiers inertes est  donc ${1\over p-1}-{1\over p(p-1)}={1\over p}$.
\bigskip

\proclaim Proposition 2. Soit $S_1$ l'ensemble des id\'eaux premiers $\qq$ de $O_K$ tels que $\qq^{{e\over p}}$ soit principal. 
Alors, $S_1$ a une densit\'e  \'egale \`a ${1\over p^N}$.
\bigskip

D\'emonstration : Un id\'eal premier $\qq$ de $O_K$ appartient \`a $S_1$ si et seulement si l'ordre de la classe de $\qq$ dans $\Cl_K$ n'est pas divisible par $p^r$. Le nombre d'\'el\'ements de $\Z/p^r\Z$ qui ne sont pas d'ordre $p^r$ est  $p^{r-1}$. Ainsi, dans 
un produit de $N$ copies de $\Z/p^r\Z$, 
il y a $p^{N(r-1)}$ \'el\'ements 
dont l'ordre n'est pas divisible par $p^r$. Le nombre d'\'el\'ements de $\Cl_K$ d'ordre non divisible par $p^r$ est donc
$$p^{N(r-1)}{h_K\over p^{rN}}={h_K\over p^N}.$$
Par ailleurs, il y une densit\'e de ${1\over h_K}$ d'id\'eaux premiers de $O_K$ dans chaque classe d'id\'eaux de $K$, d'o\`u le r\'esultat.
\bigskip

{\bf{Fin de la d\'emonstration du th\'eor\`eme.}} D'apr\`es la proposition 1, utilis\'ee avec $A=S$,  et l'hypoth\`ese faite sur $S$, l'ensemble $T$ des id\'eaux premiers de $O_K$ de degr\'e $1$ au-dessus d'un nombre premier de $S$,  est  de densit\'e strictement plus grande que ${1\over p^N}$. D'apr\`es la proposition 2, 
il existe  donc  une infinit\'e 
d'id\'eaux premiers $\qq$ de $T$ tels que $\qq^{{e\over p}}$ ne soit pas principal. Cela entra\^\i ne que $S_0$ est infini.
\vskip0pt\noindent
Il reste \`a remarquer que  si $\ell$ et $\ell'$ sont deux nombres premiers distincts,  ne divisant pas $c$, la jacobienne de $C_{\ell}$ a bonne r\'eduction en $\ell'$ mais pas en $\ell$,  et de m\^eme pour $C_{\ell'}$. En particulier, les courbes $C_{\ell}$ et $C_{\ell'}$ ne sont pas $\Q$-isomorphes, d'o\`u le th\'eor\`eme.
\bigskip

{\douze {III.  D\'emonstration du corollaire}}
\medskip
Soit $p$ un nombre premier fix\'e. Suivant la terminologie adopt\'ee dans le paragraphe 3 de [H-K], on dira qu'un nombre premier $q$ est  exceptionnel pour $p$,  si  on a $q\neq p$ et s'il existe 
$a,b,c$ non nuls dans $\FF_q$ tels que la courbe de Fermat   d'\'equation $ax^p+by^p+cz^p=0$ ne poss\`ede pas de points rationnels sur $\FF_q$.   
D'apr\`es les travaux de Weil, 
les  nombres premiers exceptionnels pour $p$ sont plus petits que $\bigl((p-1)(p-2)\bigr)^2$ ([We]). En particulier, ils sont  explicitables. Par ailleurs, ils sont congrus \`a $1$
 modulo $p$. 
Afin d'\'etablir qu'une courbe  $C_{\ell}$ n'a pas d'obstructions locales, il suffit de v\'erifier que pour tout nombre premier $q$ divisant $p\ell c$  ou qui est exceptionnel pour $p$, l'ensemble $C_{\ell}(\Q_q)$ est non vide. On a utilis\'e pour cela le lemme 3.1 de  [H-K].
\smallskip
Pour chaque couple $(p,c)$,  notons  $d$ la densit\'e de l'ensemble $S$ correspondant.  
\medskip
1) Supposons $(p,c)=(3,921)$. On a $c^2\not\equiv 1 \mod 9$  et le groupe des classes du corps  $K=\Q\bigl(\root 3 \of {921}\bigr)$ est isomorphe \`a $\Z/3\Z\times \Z/3\Z$. On a ainsi $N=2$.    Il n'y a pas de nombres premiers exceptionnels pour $p=3$. Par suite, 
 $C_{\ell}$ n'a pas d'obstructions locales si et seulement si $C_{\ell}(\Q_3)$, $C_{\ell}(\Q_{\ell})$  et $C_{\ell}(\Q_{307})$ sont non vides. 
Soit $\ell$ un nombre premier   congru \`a $2$ modulo $3$.   On v\'erifie que $C_{\ell}(\Q_3)$ est non vide si et seulement si on a
 $$\ell\equiv 2,5,8\mod 9.$$
 L'ensemble $C_{\ell}(\Q_{\ell})$ est  non vide, car tout \'el\'ement de $\FF_{\ell}$ est un cube dans $\FF_{\ell}$.
Par ailleurs, on trouve 102   classes  de nombres premiers $\ell$ modulo 307 pour lesquelles $C_{\ell}(\Q_{307})$ est  non vide.
On obtient ainsi 306  classes  de nombres premiers $\ell$ modulo 2763 pour lesquelles $C_{\ell}$ n'a pas d'obstructions locales. 
En notant $\varphi$ la fonction indicatrice d'Euler, on en d\'eduit  que l'on a 
$$d={306\over \varphi(2763)}={1\over 6}.$$
On a $d>{1\over 9}$, d'o\`u le r\'esultat dans ce cas.
\smallskip
2) Supposons $(p,c)=(5,19)$. On a $c^4\not\equiv  1 \mod 25$ et $h_K=5$, en particulier $N=1$.
Soit $\ell$ un nombre premier non congru \`a $1$ modulo $5$ et distinct de $5$. Parce que 11 est le seul nombre premier exceptionnel pour $p=5$, 
la courbe 
 $C_{\ell}$ n'a pas d'obstructions locales si et seulement si  
 $$C_{\ell}(\Q_5), C_{\ell}(\Q_{\ell}), C_{\ell}(\Q_{19}), C_{\ell}(\Q_{11})$$
 ne sont pas vides. Tel est le cas de $C_{\ell}(\Q_{\ell})$  et de $C_{\ell}(\Q_{19})$ car tout \'el\'ement de $\FF_{\ell}$ et de $\FF_{19}$ est une puissance $5$-i\`eme. 
 On constate que   $C_{\ell}(\Q_5)$ est non vide si et seulement si on a 
 $$\ell\equiv 7,8,9,12,13,17,18,19,24\mod 25.$$
Par ailleurs, $C_{\ell}(\Q_{11})$  est non vide si et seulement si on a 
 $$\ell\equiv 1,2,3,4,7,8,9,10\mod 11.$$
On obtient ainsi  72   classes  de nombres premiers $\ell$ modulo 275 pour lesquelles $C_{\ell}/\Q$ n'a pas d'obstructions locales. Il en r\'esulte que l'on a 
$$d={9\over 25}.$$ 
On a $d>{1\over 5}$, d'o\`u le r\'esultat.
\medskip
Pour  chaque couple $(p,c)$ consid\'er\'e ci-dessous, on pr\'esente les r\'esultats num\'e\-riques   dans des tableaux  \`a deux entr\'ees $q$ et $N_q$, qui se lisent de la faon suivante.
L'entier  $q$ parcourt la r\'eunion de $\big\lbrace p\big\rbrace$ et de l'ensemble des nombres premiers exceptionnels pour $p$. 
 L'entier $N_p$ est le nombre de classes de nombres premiers $\ell$ modulo $p^2$, avec $\ell\not\equiv 1 \mod p$,  pour lesquelles $C_{\ell}(\Q_p)$ est non vide. Si  $q\neq p$,  $N_q$ est le nombre de classes 
de nombres premiers $\ell$ modulo $q$  pour lesquelles $C_{\ell}(\Q_q)$ est non vide.
\vskip0pt\noindent
Par ailleurs,   $c$ est premier, on a $c\not\equiv 1 \mod p$ et  $c^{p-1}\not\equiv 1 \mod p^2$.  Pour tout  nombre premier $\ell\not\equiv 1 \mod p$ et distinct de $p$, $C_{\ell}(\Q_{\ell})$ est non vide. De m\^eme,  $C_{\ell}(\Q_{c})$ est  non vide. 
\medskip
On obtient les r\'esultats suivants.
\medskip

3) Pour  $(p,c)=(7,13)$ : on a $h_K=7$, 
\bigskip 
\centerline{\vbox{\offinterlineskip
\halign{\vrule height12pt\ \hfil#\hfil\ \vrule&&\ \hfil#\hfil\ \vrule\cr
\noalign{\hrule}
\omit\vrule height2pt\hfil\vrule&&&&\cr
$q$ &$7$ &$29$&$43$&$71$\cr
\omit\vrule height2pt\hfil\vrule&&&&\cr
\noalign{\hrule}
$N_q$ &$25$ &$16$& $30$ &$60$\cr
\omit\vrule height2pt\hfil\vrule&&&&\cr
\noalign{\hrule}}
}}
\medskip
\vskip0pt\noindent
d'o\`u   
$$d={25\times 16\times 30\times 60\over \varphi(49\times 29\times 43\times 71)}={500\over 2401}>{1\over 7}.$$
\smallskip
4) Pour  $(p,c)=(11,373)$ : $\Cl_K$ est isomorphe \`a $\Z/11\Z\times \Z/11\Z$. 
 On obtient  
\bigskip
\centerline{\vbox{\offinterlineskip
\halign{\vrule height12pt\ \hfil#\hfil\ \vrule&&\ \hfil#\hfil\ \vrule\cr
\noalign{\hrule}
\omit\vrule height2pt\hfil\vrule&&&&&&\cr
$q$ &$11$ &$23$&$67$&$89$&$199$&$419$\cr
\omit\vrule height2pt\hfil\vrule&&&&&&\cr
\noalign{\hrule}
$N_q$ &$54$ &$8$& $66$ &$56$&$162$&$380$\cr
\omit\vrule height2pt\hfil\vrule&&&&&&\cr
\noalign{\hrule}}
}}
\medskip
\vskip0pt\noindent
d'o\`u   
 $$d={13608\over 161051}>{1\over 121}.$$
\medskip

 5) Pour  $(p,c)=(13,103)$ : on a  $h_K=13$,   
\bigskip
\centerline{\vbox{\offinterlineskip
\halign{\vrule height12pt\ \hfil#\hfil\ \vrule&&\ \hfil#\hfil\ \vrule\cr
\noalign{\hrule}
\omit\vrule height2pt\hfil\vrule&&&&&&&\cr
$q$ &$13$ &$53$&$79$&$131$&$157$&$313$&$547$\cr
\omit\vrule height2pt\hfil\vrule&&&&&&&\cr
\noalign{\hrule}
$N_q$ &$99$ &$24$& $78$ &$80$&$108$&$240$&$546$\cr
\omit\vrule height2pt\hfil\vrule&&&&&&&\cr
\noalign{\hrule}}
}}
\medskip
\vskip0pt\noindent
d'o\`u  
$$d={ 35640\over 371293}>{1\over 13}.$$
\medskip

6) Pour $(p,c)=(17,1087)$ : on a  $h_K=17$, 
\bigskip
\centerline{\vbox{\offinterlineskip
\halign{\vrule height12pt\ \hfil#\hfil\ \vrule&&\ \hfil#\hfil\ \vrule\cr
\noalign{\hrule}
\omit\vrule height2pt\hfil\vrule&&&&&&&&&&&&&&\cr
$q$ &$17$ &$103$&$137$&$239$&$307$&$409$&$443$&$613$&$647$&$919$&$953$&$1021$&$1123$&$1429$\cr
\omit\vrule height2pt\hfil\vrule&&&&&&&&&&&&&&\cr
\noalign{\hrule}
$N_q$ &$165$ &$102$& $56$ &$168$&$252$&$360$&$416$&$612$&$608$&$918$&$952$&$1020$&$1056$&$1428$\cr
\omit\vrule height2pt\hfil\vrule&&&&&&&&&&&&&&\cr
\noalign{\hrule}}
}}
\medskip
\vskip0pt\noindent
d'o\`u  
$$d={745113600\over 6975757441}>{1\over 17}. $$
\medskip
7) Pour  $(p,c)=(19,37)$ : on a $h_K=19$, 
\bigskip
\centerline{\vbox{\offinterlineskip
\halign{\vrule height12pt\ \hfil#\hfil\ \vrule&&\ \hfil#\hfil\ \vrule\cr
\noalign{\hrule}
\omit\vrule height2pt\hfil\vrule&&&&&&&&&&\cr
$q$ &$19$ &$191$&$229$&$419$&$457$&$571$&$647$&$761$&$1103$&$1597$\cr
\omit\vrule height2pt\hfil\vrule&&&&&&&&&&\cr
\noalign{\hrule}
$N_q$ &$187$ &$100$& $144$ &$374$&$312$&$450$&$646$&$720$&$986$&1512\cr
\omit\vrule height2pt\hfil\vrule&&&&&&&&&&\cr
\noalign{\hrule}}
}}
\medskip
\vskip0pt\noindent
d'o\`u   
$$d={22762911600\over 322687697779}>{1\over 19}.$$
\bigskip
 Cela termine la d\'emonstration du  corollaire.
\vfill\eject

\centerline {\douze{Bibliographie}}
\bigskip
\vskip0pt\noindent
[Co-1] H. Cohen,   A Course in Computational Algebraic Number Theory, Springer  GTM {\bf{138}}, 1993.
\smallskip
\vskip0pt\noindent
[Co-2] H. Cohen,   Number Theory Volume I : Tools and Diophantine Equations, Springer  GTM {\bf{239}}, 2007.
\smallskip
\vskip0pt\noindent
[H-K] E. Halberstadt et A. Kraus, Courbes de Fermat : r\'esultats et probl\`emes,  {\it{J. reine angew. Math.}} {\bf{548}} (2002), 167-234.
\smallskip
\vskip0pt\noindent
[Ne] J. Neukirch, Class field Theory, Grundlehren der mathematischen Wissenschaften, Springer-Verlag  {\bf{280}}, 1986.
\smallskip
\vskip0pt\noindent
[Pa]  C. Batut, D. Bernardi, K. Belabas, H. Cohen et M. Olivier, PARI-GP, version 2.7.3,  
Universit\'e de Bordeaux I, (2015).
\smallskip
\vskip0pt\noindent
[Se] E. S. Selmer, The Diophantine equation $ax^3+by^3+cz^3=0$, {\it{Acta Math.}} {\bf{85}} (1951),   203-362.
\smallskip
\vskip0pt\noindent
[We] A. Weil, Sur les courbes alg\'ebriques et les vari\'et\'es qui s'en d\'eduisent, Hermann, 1948.
\bigskip
\bigskip
\line {\hfill{28 janvier 2016}}

\item{} Alain Kraus
\item{}Universit\'e de Paris VI, 
\item{} Institut de Math\'ematiques de Jussieu
\item{} 4 Place Jussieu, 75005 Paris,  
\item{} France
\medskip
\vskip0pt\noindent
\item{}e-mail : alain.kraus@imj-prg.fr 

\bigskip

\bye